\documentclass{amsart}
\usepackage{amsfonts}

\newtheorem{theorem}{Theorem}
\theoremstyle{plain}

\newtheorem{corollary}{Corollary}

\newtheorem{proposition}{Proposition}
\newtheorem{remark}{Remark}

\numberwithin{equation}{section}
\input{tcilatex}

\begin{document}
\title[Reverse of the Schwarz Inequality]{Reverses of the Schwarz Inequality
in Inner Product Spaces Generalising a Klamkin-McLenaghan Result}
\author{Sever S. Dragomir}
\address{School of Computer Science and Mathematics\\
Victoria University\\
PO Box 14428, Melbourne City\\
VIC 8001, Australia.}
\email{sever.dragomir@vu.edu.au}
\urladdr{http://rgmia.vu.edu.au/dragomir}
\date{July 27, 2005}

\begin{abstract}
New reverses of the Schwarz inequality in inner product spaces that
incorporate the classical Klamkin-McLenaghan result for the case of positive 
$n-$tuples are given. Applications for Lebesgue integrals are also provided.
\end{abstract}

\maketitle

\section{Introduction}

In 2004, the author \cite{SSD1} (see also \cite{SSD2}) proved the following
reverse of the Schwarz inequality:

\begin{theorem}
\label{ta}Let $\left( H;\left\langle \cdot ,\cdot \right\rangle \right) $ be
an inner product space over the real or complex number field $\mathbb{K}$
and $x,a\in H,$ $r>0$ such that%
\begin{equation}
\left\Vert x-a\right\Vert \leq r<\left\Vert a\right\Vert .  \label{1.1}
\end{equation}%
Then%
\begin{equation}
\left\Vert x\right\Vert \left( \left\Vert a\right\Vert ^{2}-r^{2}\right) ^{%
\frac{1}{2}}\leq \func{Re}\left\langle x,a\right\rangle  \label{1.2}
\end{equation}%
or, equivalently,%
\begin{equation}
\left\Vert x\right\Vert ^{2}\left\Vert a\right\Vert ^{2}-\left[ \func{Re}%
\left\langle x,a\right\rangle \right] ^{2}\leq r^{2}\left\Vert x\right\Vert
^{2}.  \label{1.3}
\end{equation}%
The case of equality holds in (\ref{1.2}) or (\ref{1.3}) if and only if%
\begin{equation}
\left\Vert x-a\right\Vert =r\quad \text{and}\quad \left\Vert x\right\Vert
^{2}+r^{2}=\left\Vert a\right\Vert ^{2}.  \label{1.4}
\end{equation}
\end{theorem}

If above one chooses%
\begin{equation*}
a=\frac{\Gamma +\gamma }{2}\cdot y\quad \text{and}\quad r=\frac{1}{2}%
\left\vert \Gamma -\gamma \right\vert \left\Vert y\right\Vert
\end{equation*}%
then the condition (\ref{1.1}) is equivalent to%
\begin{equation}
\left\Vert x-\frac{\Gamma +\gamma }{2}\cdot y\right\Vert \leq \frac{1}{2}%
\left\vert \Gamma -\gamma \right\vert \left\Vert y\right\Vert \quad \text{and%
}\quad \func{Re}\left( \Gamma \bar{\gamma}\right) >0.  \label{1.5}
\end{equation}%
Therefore, we can state the following particular result as well:

\begin{corollary}
\label{cb}Let $\left( H;\left\langle \cdot ,\cdot \right\rangle \right) $ be
as above, $x,y\in H$ and $\gamma ,\Gamma \in \mathbb{K}$ with $\func{Re}%
\left( \Gamma \bar{\gamma}\right) >0.$ If%
\begin{equation}
\left\Vert x-\frac{\Gamma +\gamma }{2}\cdot y\right\Vert \leq \frac{1}{2}%
\left\vert \Gamma -\gamma \right\vert \left\Vert y\right\Vert  \label{1.6}
\end{equation}%
or, equivalently,%
\begin{equation}
\func{Re}\left\langle \Gamma y-x,x-\gamma y\right\rangle \geq 0,  \label{1.7}
\end{equation}%
then%
\begin{align}
\left\Vert x\right\Vert \left\Vert y\right\Vert & \leq \frac{\func{Re}\left[
\left( \bar{\Gamma}+\bar{\gamma}\right) \left\langle x,a\right\rangle \right]
}{2\sqrt{\func{Re}\left( \Gamma \bar{\gamma}\right) }}  \label{1.8} \\
& =\frac{\func{Re}\left( \Gamma +\gamma \right) \func{Re}\left\langle
x,a\right\rangle +\func{Im}\left( \Gamma +\gamma \right) \func{Im}%
\left\langle x,a\right\rangle }{2\sqrt{\func{Re}\left( \Gamma \bar{\gamma}%
\right) }}  \notag \\
& \left( \leq \frac{\left\vert \Gamma +\gamma \right\vert }{\sqrt{\func{Re}%
\left( \Gamma \bar{\gamma}\right) }}\left\vert \left\langle x,a\right\rangle
\right\vert \right) .  \notag
\end{align}%
The case of equality holds in (\ref{1.8}) if and only if the equality case
holds in (\ref{1.6}) (or (\ref{1.7})) and%
\begin{equation}
\left\Vert x\right\Vert =\sqrt{\func{Re}\left( \Gamma \bar{\gamma}\right) }%
\left\Vert y\right\Vert .  \label{1.9}
\end{equation}
\end{corollary}

If the restriction $\left\Vert a\right\Vert >r$ is removed from Theorem \ref%
{ta}, then a different reverse of the Schwarz inequality may be stated \cite%
{SSD3} (see also \cite{SSD2}):

\begin{theorem}
\label{tc}Let $\left( H;\left\langle \cdot ,\cdot \right\rangle \right) $ be
an inner product space over $\mathbb{K}$ and $x,a\in H,$ $r>0$ such that%
\begin{equation}
\left\Vert x-a\right\Vert \leq r.  \label{1.10}
\end{equation}%
Then%
\begin{equation}
\left\Vert x\right\Vert \left\Vert a\right\Vert -\func{Re}\left\langle
x,a\right\rangle \leq \frac{1}{2}r^{2}.  \label{1.11}
\end{equation}%
The equality holds in (\ref{1.11}) if and only if the equality case is
realised in (\ref{1.10}) and $\left\Vert x\right\Vert =\left\Vert
a\right\Vert .$
\end{theorem}

As a corollary of the above, we can state:

\begin{corollary}
\label{cd}Let $\left( H;\left\langle \cdot ,\cdot \right\rangle \right) $ be
as above, $x,y\in H$ and $\gamma ,\Gamma \in \mathbb{K}$ with $\Gamma \neq
-\gamma .$ If either (\ref{1.6}) or, equivalently, (\ref{1.7}) hold true,
then%
\begin{equation}
\left\Vert x\right\Vert \left\Vert y\right\Vert -\frac{\func{Re}\left(
\Gamma +\gamma \right) \func{Re}\left\langle x,a\right\rangle +\func{Im}%
\left( \Gamma +\gamma \right) \func{Im}\left\langle x,a\right\rangle }{%
\left\vert \Gamma +\gamma \right\vert }\leq \frac{1}{4}\cdot \frac{%
\left\vert \Gamma -\gamma \right\vert ^{2}}{\left\vert \Gamma +\gamma
\right\vert }\left\Vert y\right\Vert ^{2}.  \label{1.12}
\end{equation}%
The equality holds in (\ref{1.12}) if and only if the equality case is
realised in either (\ref{1.6}) or (\ref{1.7}) and%
\begin{equation}
\left\Vert x\right\Vert =\frac{1}{2}\left\vert \Gamma +\gamma \right\vert
\left\Vert y\right\Vert .  \label{1.13}
\end{equation}
\end{corollary}

As pointed out in \cite{SSD4}, the above results are motivated by the fact
that they generalise to the case of real or complex inner product spaces
some classical reverses of the Cauchy-Bunyakovsky-Schwarz inequality for
positive $n-$tuples due to Polya-Szeg\"{o} \cite{PS}, Cassels \cite{W},
Shisha-Mond \cite{SM} and Greub-Rheinboldt \cite{GR}.

The main aim of this paper is to establish a new reverse of Schwarz's
inequality similar to the ones in Theorems \ref{ta} and \ref{tc} which will
reduce, for the particular case of positive $n-$tuples, to the
Klamkin-McLenaghan result from \cite{KM}.

\section{The Results}

The following result may be stated.

\begin{theorem}
\label{t2.1}Let $\left( H;\left\langle \cdot ,\cdot \right\rangle \right) $
be an inner product space over the real or complex number field $\mathbb{K}$
and $x,a\in H,$ $r>0$ with $\left\langle x,a\right\rangle \neq 0$ and%
\begin{equation}
\left\Vert x-a\right\Vert \leq r<\left\Vert a\right\Vert .  \label{2.1}
\end{equation}%
Then%
\begin{equation}
\frac{\left\Vert x\right\Vert ^{2}}{\left\vert \left\langle x,a\right\rangle
\right\vert }-\frac{\left\vert \left\langle x,a\right\rangle \right\vert }{%
\left\Vert a\right\Vert ^{2}}\leq \frac{2r^{2}}{\left\Vert a\right\Vert
\left( \left\Vert a\right\Vert +\sqrt{\left\Vert a\right\Vert ^{2}-r^{2}}%
\right) },  \label{2.2}
\end{equation}%
with equality if and only if the equality case holds in (\ref{2.1}) and%
\begin{equation}
\func{Re}\left\langle x,a\right\rangle =\left\vert \left\langle
x,a\right\rangle \right\vert =\left\Vert a\right\Vert \left( \left\Vert
a\right\Vert ^{2}-r^{2}\right) ^{\frac{1}{2}}.  \label{2.3}
\end{equation}%
The constant $2$ is best possible in (\ref{2.2}) in the sese that it cannot
be replaced by a smaller quantity.
\end{theorem}

\begin{proof}
The first condition in (\ref{2.1}) is obviously equivalent with%
\begin{equation}
\frac{\left\Vert x\right\Vert ^{2}}{\left\vert \left\langle x,a\right\rangle
\right\vert }\leq \frac{2\func{Re}\left\langle x,a\right\rangle }{\left\vert
\left\langle x,a\right\rangle \right\vert }-\frac{\left\Vert a\right\Vert
^{2}-r^{2}}{\left\vert \left\langle x,a\right\rangle \right\vert }
\label{2.4}
\end{equation}%
with equality if and only if $\left\Vert x-a\right\Vert =r.$

Subtracting from both sides of (\ref{2.4}) the same quantity $\frac{%
\left\vert \left\langle x,a\right\rangle \right\vert }{\left\Vert
a\right\Vert ^{2}}$ and performing some elementary calculations, we get the
equivalent inequality:%
\begin{multline}
\frac{\left\Vert x\right\Vert ^{2}}{\left\vert \left\langle x,a\right\rangle
\right\vert }-\frac{\left\vert \left\langle x,a\right\rangle \right\vert }{%
\left\Vert a\right\Vert ^{2}}  \label{2.5} \\
\leq 2\cdot \frac{\func{Re}\left\langle x,a\right\rangle }{\left\vert
\left\langle x,a\right\rangle \right\vert }-\left( \frac{\left\vert
\left\langle x,a\right\rangle \right\vert ^{\frac{1}{2}}}{\left\Vert
a\right\Vert }-\frac{\left( \left\Vert a\right\Vert ^{2}-r^{2}\right) ^{%
\frac{1}{2}}}{\left\vert \left\langle x,a\right\rangle \right\vert ^{\frac{1%
}{2}}}\right) ^{2}-\frac{2\sqrt{\left\Vert a\right\Vert ^{2}-r^{2}}}{%
\left\Vert a\right\Vert }.
\end{multline}%
Since, obviously%
\begin{equation*}
\func{Re}\left\langle x,a\right\rangle \leq \left\vert \left\langle
x,a\right\rangle \right\vert \quad \text{and}\quad \left( \frac{\left\vert
\left\langle x,a\right\rangle \right\vert ^{\frac{1}{2}}}{\left\Vert
a\right\Vert }-\frac{\left( \left\Vert a\right\Vert ^{2}-r^{2}\right) ^{%
\frac{1}{2}}}{\left\vert \left\langle x,a\right\rangle \right\vert ^{\frac{1%
}{2}}}\right) ^{2}\geq 0,
\end{equation*}%
hence, by (\ref{2.5}) we get%
\begin{equation}
\frac{\left\Vert x\right\Vert ^{2}}{\left\vert \left\langle x,a\right\rangle
\right\vert }-\frac{\left\vert \left\langle x,a\right\rangle \right\vert }{%
\left\Vert a\right\Vert ^{2}}\leq 2\left( 1-\frac{\sqrt{\left\Vert
a\right\Vert ^{2}-r^{2}}}{\left\Vert a\right\Vert }\right)  \label{2.6}
\end{equation}%
with equality if and only if%
\begin{equation}
\left\Vert x-a\right\Vert =r,\quad \func{Re}\left\langle x,a\right\rangle
=\left\vert \left\langle x,a\right\rangle \right\vert \quad \text{and}\quad
\left\vert \left\langle x,a\right\rangle \right\vert =\left\Vert
a\right\Vert \left( \left\Vert a\right\Vert ^{2}-r^{2}\right) ^{\frac{1}{2}}.
\label{2.7}
\end{equation}%
Observe that (\ref{2.6}) is equivalent with (\ref{2.2}) and the first part
of the theorem is proved.

To prove the sharpness of the constant, let us assume that there is a $C>0$
such that%
\begin{equation}
\frac{\left\Vert x\right\Vert ^{2}}{\left\vert \left\langle x,a\right\rangle
\right\vert }-\frac{\left\vert \left\langle x,a\right\rangle \right\vert }{%
\left\Vert a\right\Vert ^{2}}\leq \frac{Cr^{2}}{\left\Vert a\right\Vert
\left( \left\Vert a\right\Vert +\sqrt{\left\Vert a\right\Vert ^{2}-r^{2}}%
\right) },  \label{2.7.1}
\end{equation}%
provided $\left\Vert x-a\right\Vert \leq r<\left\Vert a\right\Vert .$

Now, consider $\varepsilon \in \left( 0,1\right) $ and let $r=\sqrt{%
\varepsilon },$ $a,e\in H,$ $\left\Vert a\right\Vert =\left\Vert
e\right\Vert =1$ and $a\perp e.$ Define $x:=a+\sqrt{\varepsilon }e.$ We
observe that $\left\Vert x-a\right\Vert =\sqrt{\varepsilon }=r<1=\left\Vert
a\right\Vert ,$ which shows that the condition (\ref{2.1}) of the theorem is
fulfilled. We also observe that%
\begin{equation*}
\left\Vert x\right\Vert ^{2}=\left\Vert a\right\Vert ^{2}+\varepsilon
\left\Vert e\right\Vert ^{2}=1+\varepsilon ,\quad \left\langle
x,a\right\rangle =\left\Vert e\right\Vert ^{2}=1
\end{equation*}%
and utilising (\ref{2.7.1}) we get%
\begin{equation*}
1+\varepsilon -1\leq \frac{C\varepsilon }{\left( 1+\sqrt{1-\varepsilon }%
\right) },
\end{equation*}%
giving $1+\sqrt{1-\varepsilon }\leq C$ for any $\varepsilon \in \left(
0,1\right) .$ Letting $\varepsilon \rightarrow 0+,$ we get $C\geq 2,$ which
shows that the constant $2$ in (\ref{2.2}) is best possible.
\end{proof}

\begin{remark}
\label{r2.1}In a similar manner, one can prove that if $\func{Re}%
\left\langle x,a\right\rangle \neq 0$ and (\ref{2.2}) holds true, then:%
\begin{equation}
\frac{\left\Vert x\right\Vert ^{2}}{\left\vert \func{Re}\left\langle
x,a\right\rangle \right\vert }-\frac{\left\vert \func{Re}\left\langle
x,a\right\rangle \right\vert }{\left\Vert a\right\Vert ^{2}}\leq \frac{2r^{2}%
}{\left\Vert a\right\Vert \left( \left\Vert a\right\Vert +\sqrt{\left\Vert
a\right\Vert ^{2}-r^{2}}\right) }  \label{2.8}
\end{equation}%
with equality if and only if $\left\Vert x-a\right\Vert =r$ and%
\begin{equation}
\func{Re}\left\langle x,a\right\rangle =\left\Vert a\right\Vert \left(
\left\Vert a\right\Vert ^{2}-r^{2}\right) ^{\frac{1}{2}}.  \label{2.9}
\end{equation}%
The constant $2$ is best possible in (\ref{2.8}).
\end{remark}

\begin{remark}
\label{r2.2}Since (\ref{2.2}) is equivalent with%
\begin{equation}
\left\Vert x\right\Vert ^{2}\left\Vert a\right\Vert ^{2}-\left\vert
\left\langle x,a\right\rangle \right\vert ^{2}\leq \frac{2r^{2}\left\Vert
a\right\Vert ^{2}}{\left\Vert a\right\Vert \left( \left\Vert a\right\Vert +%
\sqrt{\left\Vert a\right\Vert ^{2}-r^{2}}\right) }\left\vert \left\langle
x,a\right\rangle \right\vert  \label{2.10}
\end{equation}%
and (\ref{2.8}) is equivalent to%
\begin{equation}
\left\Vert x\right\Vert ^{2}\left\Vert a\right\Vert ^{2}-\left[ \func{Re}%
\left\langle x,a\right\rangle \right] ^{2}\leq \frac{2r^{2}\left\Vert
a\right\Vert ^{2}}{\left\Vert a\right\Vert \left( \left\Vert a\right\Vert +%
\sqrt{\left\Vert a\right\Vert ^{2}-r^{2}}\right) }\left\vert \func{Re}%
\left\langle x,a\right\rangle \right\vert  \label{2.11}
\end{equation}%
hence (\ref{2.11}) is a tighter inequality than (\ref{2.10}), because in
complex spaces, in general $\left\vert \left\langle x,a\right\rangle
\right\vert >\left\vert \func{Re}\left\langle x,a\right\rangle \right\vert .$
\end{remark}

The following corollary is of interest.

\begin{corollary}
\label{c2.1}Let $\left( H;\left\langle \cdot ,\cdot \right\rangle \right) $
be a real or complex inner product space and $x,y\in H$ with $\left\langle
x,y\right\rangle \neq 0$, $\gamma ,\Gamma \in \mathbb{K}$ with $\func{Re}%
\left( \Gamma \bar{\gamma}\right) >0.$ If either (\ref{2.6}) or,
equivalently (\ref{2.7}) holds true, then%
\begin{equation}
\frac{\left\Vert x\right\Vert ^{2}}{\left\vert \left\langle x,y\right\rangle
\right\vert }-\frac{\left\vert \left\langle x,y\right\rangle \right\vert }{%
\left\Vert y\right\Vert ^{2}}\leq \left\vert \Gamma +\gamma \right\vert -2%
\sqrt{\func{Re}\left( \Gamma \bar{\gamma}\right) }.  \label{2.12}
\end{equation}%
The equality holds in (\ref{2.12}) if and only if the equality case holds in
(\ref{2.6}) (or in (\ref{2.7})) and%
\begin{equation}
\func{Re}\left[ \left( \Gamma +\gamma \right) \left\langle x,y\right\rangle %
\right] =\left\vert \Gamma +\gamma \right\vert \left\vert \left\langle
x,y\right\rangle \right\vert =\left\vert \Gamma +\gamma \right\vert \sqrt{%
\func{Re}\left( \Gamma \bar{\gamma}\right) }\left\Vert y\right\Vert ^{2}.
\label{2.13}
\end{equation}
\end{corollary}

\begin{proof}
We use the inequality (\ref{2.2}) in its equivalent form%
\begin{equation*}
\frac{\left\Vert x\right\Vert ^{2}}{\left\vert \left\langle x,a\right\rangle
\right\vert }-\frac{\left\vert \left\langle x,a\right\rangle \right\vert }{%
\left\Vert a\right\Vert ^{2}}\leq \frac{2\left( \left\Vert a\right\Vert +%
\sqrt{\left\Vert a\right\Vert ^{2}-r^{2}}\right) }{\left\Vert a\right\Vert }.
\end{equation*}%
Choosing $a=\frac{\Gamma +\gamma }{2}\cdot y$ and $r=\frac{1}{2}\left\vert
\Gamma -\gamma \right\vert \left\Vert y\right\Vert ,$ we have%
\begin{multline*}
\frac{\left\Vert x\right\Vert ^{2}}{\left\vert \frac{\Gamma +\gamma }{2}%
\right\vert \left\vert \left\langle x,y\right\rangle \right\vert }-\frac{%
\left\vert \frac{\Gamma +\gamma }{2}\right\vert \left\vert \left\langle
x,y\right\rangle \right\vert }{\left\vert \frac{\Gamma +\gamma }{2}%
\right\vert ^{2}\left\Vert y\right\Vert ^{2}} \\
\leq \frac{2\left( \left\vert \frac{\Gamma +\gamma }{2}\right\vert
\left\vert \left\Vert y\right\Vert \right\vert +\sqrt{\left\vert \frac{%
\Gamma +\gamma }{2}\right\vert ^{2}\left\Vert y\right\Vert ^{2}-\frac{1}{4}%
\left\vert \Gamma -\gamma \right\vert ^{2}\left\Vert y\right\Vert ^{2}}%
\right) }{\left\vert \frac{\Gamma -\gamma }{2}\right\vert \left\Vert
y\right\Vert }
\end{multline*}%
which is equivalent to (\ref{2.12}).
\end{proof}

\begin{remark}
\label{r2.3}The inequality (\ref{2.12}) has been obtained in a different way
in \cite[Theorem 2]{EMP}. However, in \cite{EMP} the authors did not
consider the equality case which may be of interest for applications.
\end{remark}

\begin{remark}
\label{r2.4}If we assume that $\Gamma =M\geq m=\gamma >0,$ which is very
convenient in applications, then%
\begin{equation}
\frac{\left\Vert x\right\Vert ^{2}}{\left\vert \left\langle x,y\right\rangle
\right\vert }-\frac{\left\vert \left\langle x,y\right\rangle \right\vert }{%
\left\Vert y\right\Vert ^{2}}\leq \left( \sqrt{M}-\sqrt{m}\right) ^{2},
\label{2.14}
\end{equation}%
provided that either%
\begin{equation}
\func{Re}\left\langle My-x,x-my\right\rangle \geq 0  \label{2.15}
\end{equation}%
or, equivalently,%
\begin{equation}
\left\Vert x-\frac{m+M}{2}y\right\Vert \leq \frac{1}{2}\left( M-m\right)
\left\Vert y\right\Vert   \label{2.16}
\end{equation}%
holds true.
\end{remark}

The equality holds in (\ref{2.14}) if and only if the equality case holds in
(\ref{2.15}) (or in (\ref{2.16})) and%
\begin{equation}
\func{Re}\left\langle x,y\right\rangle =\left\vert \left\langle
x,y\right\rangle \right\vert =\sqrt{Mm}\left\Vert y\right\Vert ^{2}.
\label{2.17}
\end{equation}%
The multiplicative constant $C=1$ in front of $\left( \sqrt{M}-\sqrt{m}%
\right) ^{2}$ cannot be replaced in general with a smaller positive quantity.

Now for a non-zero complex number $z,$ we define $\func{sgn}\left( z\right)
:=\frac{z}{\left\vert z\right\vert }.$

The following result may be stated:

\begin{proposition}
\label{p2.1}Let $\left( H;\left\langle \cdot ,\cdot \right\rangle \right) $
be a real or complex inner product space and $x,y\in H$ with $\func{Re}%
\left\langle x,y\right\rangle \neq 0$ and $\gamma ,\Gamma \in \mathbb{K}$
with $\func{Re}\left( \Gamma \bar{\gamma}\right) >0.$ If either (\ref{2.6})
or, equivalently, (\ref{2.7}) hold true, then%
\begin{align}
(0& \leq \left\Vert x\right\Vert ^{2}\left\Vert y\right\Vert ^{2}-\left\vert
\left\langle x,y\right\rangle \right\vert ^{2}\leq )  \label{2.18} \\
& \left\Vert x\right\Vert ^{2}\left\Vert y\right\Vert ^{2}-\left[ \func{Re}%
\left( \func{sgn}\left( \frac{\Gamma +\gamma }{2}\right) \cdot \left\langle
x,y\right\rangle \right) \right] ^{2}  \notag \\
& \leq \left( \left\vert \Gamma +\gamma \right\vert -2\sqrt{\func{Re}\left(
\Gamma \bar{\gamma}\right) }\right) \left\vert \func{Re}\left( \func{sgn}%
\left( \frac{\Gamma +\gamma }{2}\right) \cdot \left\langle x,y\right\rangle
\right) \right\vert \left\Vert y\right\Vert ^{2}  \notag \\
& \left( \leq \left( \left\vert \Gamma +\gamma \right\vert -2\sqrt{\func{Re}%
\left( \Gamma \bar{\gamma}\right) }\right) \left\vert \left\langle
x,y\right\rangle \right\vert \left\Vert y\right\Vert ^{2}\right) .  \notag
\end{align}%
The equality holds in (\ref{2.18}) if and only if the equality case holds in
(\ref{2.6}) (or in (\ref{2.7})) and%
\begin{equation*}
\func{Re}\left[ \func{sgn}\left( \frac{\Gamma +\gamma }{2}\right) \cdot
\left\langle x,y\right\rangle \right] =\sqrt{\func{Re}\left( \Gamma \bar{%
\gamma}\right) }\left\Vert y\right\Vert ^{2}.
\end{equation*}
\end{proposition}

\begin{proof}
The inequality (\ref{2.8}) is equivalent with:%
\begin{equation*}
\left\Vert x\right\Vert ^{2}\left\Vert a\right\Vert ^{2}-\left[ \func{Re}%
\left\langle x,a\right\rangle \right] ^{2}\leq 2\left( \left\Vert
a\right\Vert +\sqrt{\left\Vert a\right\Vert ^{2}-r^{2}}\right) \cdot
\left\vert \func{Re}\left\langle x,a\right\rangle \right\vert \left\Vert
a\right\Vert .
\end{equation*}%
If in this inequality we choose $a=\frac{\Gamma +\gamma }{2}\cdot y$ \ and $%
r=\frac{1}{2}\left\vert \Gamma -\gamma \right\vert |y$ , we have%
\begin{multline*}
\left\Vert x\right\Vert ^{2}\left\vert \frac{\Gamma +\gamma }{2}\right\vert
^{2}\left\Vert y\right\Vert ^{2}-\left( \func{Re}\left[ \left( \frac{\Gamma
+\gamma }{2}\right) \cdot \left\langle x,y\right\rangle \right] \right) ^{2}
\\
\leq 2\left( \left\vert \frac{\Gamma +\gamma }{2}\right\vert \left\Vert
y\right\Vert -\sqrt{\left\vert \frac{\Gamma +\gamma }{2}\right\vert
^{2}\left\Vert y\right\Vert ^{2}-\frac{1}{4}\left\vert \Gamma -\gamma
\right\vert ^{2}\left\Vert y\right\Vert ^{2}}\right) \\
\times \left\vert \func{Re}\left[ \left( \frac{\Gamma +\gamma }{2}\right)
\cdot \left\langle x,y\right\rangle \right] \right\vert \left\vert \frac{%
\Gamma +\gamma }{2}\right\vert \left\Vert y\right\Vert ,
\end{multline*}%
which, on dividing by $\left\vert \frac{\Gamma +\gamma }{2}\right\vert
^{2}\neq 0$ (since $\func{Re}\left( \Gamma \bar{\gamma}\right) >0$), is
clearly equivalent to (\ref{2.18}).
\end{proof}

\begin{remark}
\label{r2.5}If we assume that $x,y,m,M$ satisfy either (\ref{2.15}) or,
equivalently (\ref{2.16}), then%
\begin{equation}
\frac{\left\Vert x\right\Vert ^{2}}{\left\vert \func{Re}\left\langle
x,y\right\rangle \right\vert }-\frac{\left\vert \func{Re}\left\langle
x,y\right\rangle \right\vert }{\left\Vert y\right\Vert ^{2}}\leq \left( 
\sqrt{M}-\sqrt{m}\right) ^{2}  \label{2.19}
\end{equation}%
or, equivalently%
\begin{equation}
\left\Vert x\right\Vert ^{2}\left\Vert y\right\Vert ^{2}-\left[ \func{Re}%
\left\langle x,y\right\rangle \right] ^{2}\leq \left( \sqrt{M}-\sqrt{m}%
\right) ^{2}\left\vert \func{Re}\left\langle x,y\right\rangle \right\vert
\left\Vert y\right\Vert ^{2}.  \label{2.20}
\end{equation}%
The equality holds in (\ref{2.19}) (or (\ref{2.20})) if and only if the case
of equality is valid in (\ref{2.15}) (or (\ref{2.16})) and%
\begin{equation}
\func{Re}\left\langle x,y\right\rangle =\sqrt{Mm}\left\Vert y\right\Vert
^{2}.  \label{2.21}
\end{equation}
\end{remark}

\section{Applications for Integrals}

Let $\left( \Omega ,\Sigma ,\mu \right) $ be a measure space consisting of a
set $\Omega ,$ a $\sigma -$algebra of parts $\Sigma $ and a countably
additive and positive measure $\mu $ on $\Sigma $ with values in $\mathbb{R}%
\cup \left\{ \infty \right\} .$

Denote by $L_{\rho }^{2}\left( \Omega ,\mathbb{K}\right) $ the Hilbert space
of all $\mathbb{K}$--valued functions $f$ defined on $\Omega $ that are $%
2-\rho -$integrable on $\Omega ,$ i.e.,$\int_{\Omega }\rho \left( t\right)
\left\vert f\left( s\right) \right\vert ^{2}d\mu \left( s\right) <\infty ,$
where $\rho :\Omega \rightarrow \lbrack 0,\infty )$ is a measurable function
on $\Omega .$

The following proposition contains a reverse of the
Cauchy-Bunyakovsky-Schwarz integral inequality:

\begin{proposition}
\label{p3.1}Let $f,g\in L_{\rho }^{2}\left( \Omega ,\mathbb{K}\right) $, $%
r>0 $ be such that%
\begin{equation}
\int_{\Omega }\rho \left( t\right) \left\vert f\left( t\right) -g\left(
t\right) \right\vert ^{2}d\mu \left( t\right) \leq r^{2}<\int_{\Omega }\rho
\left( t\right) \left\vert g\left( t\right) \right\vert ^{2}d\mu \left(
t\right) .  \label{3.0}
\end{equation}%
Then%
\begin{multline}
\int_{\Omega }\rho \left( t\right) \left\vert f\left( t\right) \right\vert
^{2}d\mu \left( t\right) \int_{\Omega }\rho \left( t\right) \left\vert
g\left( t\right) \right\vert ^{2}d\mu \left( t\right) -\left\vert
\int_{\Omega }\rho \left( t\right) f\left( t\right) \overline{g\left(
t\right) }d\mu \left( t\right) \right\vert ^{2}  \label{3.1} \\
\leq 2\left( \int_{\Omega }\rho \left( t\right) \left\vert g\left( t\right)
\right\vert ^{2}d\mu \left( t\right) \right) ^{\frac{1}{2}}\left\vert
\int_{\Omega }\rho \left( t\right) f\left( t\right) \overline{g\left(
t\right) }d\mu \left( t\right) \right\vert \\
\times \left[ \left( \int_{\Omega }\rho \left( t\right) \left\vert g\left(
t\right) \right\vert ^{2}d\mu \left( t\right) \right) ^{\frac{1}{2}}-\left(
\int_{\Omega }\rho \left( t\right) \left\vert g\left( t\right) \right\vert
^{2}d\mu \left( t\right) -r^{2}\right) ^{\frac{1}{2}}\right] .
\end{multline}%
The constant $2$ is sharp in (\ref{3.1}).
\end{proposition}

The proof follows from Theorem \ref{t2.1} applied for the Hilbert space $%
\left( L_{\rho }^{2}\left( \Omega ,\mathbb{K}\right) ,\left\langle \cdot
,\cdot \right\rangle _{\rho }\right) $ where%
\begin{equation*}
\left\langle f,g\right\rangle _{\rho }:=\int_{\Omega }\rho \left( t\right)
f\left( t\right) \overline{g\left( t\right) }d\mu \left( t\right) .
\end{equation*}

\begin{remark}
\label{r3.1}We observe that if $\int_{\Omega }\rho \left( t\right) d\mu
\left( t\right) =1,$ then a simple sufficient condition for (\ref{3.0}) to
hold is%
\begin{equation}
\left\vert f\left( t\right) -g\left( t\right) \right\vert \leq r<\left\vert
g\left( t\right) \right\vert \quad \text{for \ }\mu -\text{a.e. \ }t\in
\Omega .  \label{3.2}
\end{equation}
\end{remark}

The second general integral inequality is incorporated in:

\begin{proposition}
\label{p3.2}Let $f,g\in L_{\rho }^{2}\left( \Omega ,\mathbb{K}\right) $ and $%
\Gamma ,\gamma \in \mathbb{K}$ with $\func{Re}\left( \Gamma \bar{\gamma}%
\right) >0.$ If either%
\begin{equation}
\int_{\Omega }\func{Re}\left[ \left( \Gamma g\left( t\right) -f\left(
t\right) \right) \left( \overline{f\left( t\right) }-\bar{\gamma}\overline{%
g\left( t\right) }\right) \right] \rho \left( t\right) d\mu \left( t\right)
\geq 0  \label{3.3}
\end{equation}%
or, equivalently,%
\begin{multline}
\left( \int_{\Omega }\rho \left( t\right) \left\vert f\left( t\right) -\frac{%
\Gamma +\gamma }{2}g\left( t\right) \right\vert ^{2}d\mu \left( t\right)
\right) ^{\frac{1}{2}}  \label{3.4} \\
\leq \frac{1}{2}\left\vert \Gamma -\gamma \right\vert \left( \int_{\Omega
}\rho \left( t\right) \left\vert g\left( t\right) \right\vert ^{2}d\mu
\left( t\right) \right) ^{\frac{1}{2}}
\end{multline}%
holds, then%
\begin{multline}
\int_{\Omega }\rho \left( t\right) \left\vert f\left( t\right) \right\vert
^{2}d\mu \left( t\right) \int_{\Omega }\rho \left( t\right) \left\vert
g\left( t\right) \right\vert ^{2}d\mu \left( t\right) -\left\vert
\int_{\Omega }\rho \left( t\right) f\left( t\right) \overline{g\left(
t\right) }d\mu \left( t\right) \right\vert ^{2}  \label{3.5} \\
\leq \left[ \left\vert \Gamma +\gamma \right\vert -2\sqrt{\func{Re}\left(
\Gamma \bar{\gamma}\right) }\right] \left\vert \int_{\Omega }\rho \left(
t\right) f\left( t\right) \overline{g\left( t\right) }d\mu \left( t\right)
\right\vert \int_{\Omega }\rho \left( t\right) \left\vert g\left( t\right)
\right\vert ^{2}d\mu \left( t\right) .
\end{multline}
\end{proposition}

The proof is obvious by Corollary \ref{c2.1}.

\begin{remark}
\label{r3.2}A simple sufficient condition for the inequality (\ref{3.3}) to
hold is:%
\begin{equation}
\func{Re}\left[ \left( \Gamma g\left( t\right) -f\left( t\right) \right)
\left( \overline{f\left( t\right) }-\bar{\gamma}\overline{g\left( t\right) }%
\right) \right] \geq 0,  \label{3.6}
\end{equation}%
for $\mu -$a.e. $t\in \Omega .$
\end{remark}

A more convenient result that may be useful in applications is:

\begin{corollary}
\label{c3.3}If $f,g\in L_{\rho }^{2}\left( \Omega ,\mathbb{K}\right) $ and $%
M\geq m>0$ such that either%
\begin{equation}
\int_{\Omega }\func{Re}\left[ \left( Mg\left( t\right) -f\left( t\right)
\right) \left( \overline{f\left( t\right) }-m\overline{g\left( t\right) }%
\right) \right] f\left( t\right) d\mu \left( t\right) \geq 0  \label{3.7}
\end{equation}%
or, equivalently,%
\begin{multline}
\left( \int_{\Omega }\rho \left( t\right) \left\vert f\left( t\right) -\frac{%
M+m}{2}g\left( t\right) \right\vert ^{2}d\mu \left( t\right) \right) ^{\frac{%
1}{2}}  \label{3.8} \\
\leq \frac{1}{2}\left( M-m\right) \left( \int_{\Omega }\rho \left( t\right)
\left\vert g\left( t\right) \right\vert ^{2}d\mu \left( t\right) \right) ^{%
\frac{1}{2}},
\end{multline}%
holds, then%
\begin{multline}
\int_{\Omega }\rho \left( t\right) \left\vert f\left( t\right) \right\vert
^{2}d\mu \left( t\right) \int_{\Omega }\rho \left( t\right) \left\vert
g\left( t\right) \right\vert ^{2}d\mu \left( t\right) -\left\vert
\int_{\Omega }\rho \left( t\right) f\left( t\right) \overline{g\left(
t\right) }d\mu \left( t\right) \right\vert ^{2}  \label{3.9} \\
\leq \left( \sqrt{M}-\sqrt{m}\right) ^{2}\left\vert \int_{\Omega }\rho
\left( t\right) f\left( t\right) \overline{g\left( t\right) }d\mu \left(
t\right) \right\vert \int_{\Omega }\rho \left( t\right) \left\vert g\left(
t\right) \right\vert ^{2}d\mu \left( t\right) .
\end{multline}
\end{corollary}

\begin{remark}
\label{r3.3}Since, obviously,%
\begin{multline*}
\func{Re}\left[ \left( Mg\left( t\right) -f\left( t\right) \right) \left( 
\overline{f\left( t\right) }-m\overline{g\left( t\right) }\right) \right] \\
=\left( M\func{Re}g\left( t\right) -\func{Re}f\left( t\right) \right) \left( 
\func{Re}f\left( t\right) -m\func{Re}g\left( t\right) \right) \\
+\left( M\func{Im}g\left( t\right) -\func{Im}f\left( t\right) \right) \left( 
\func{Im}f\left( t\right) -m\func{Im}g\left( t\right) \right)
\end{multline*}%
for any $t\in \Omega ,$ hence a very simple sufficient condition that can be
useful in practical applications for (\ref{3.7}) to hold is:%
\begin{equation*}
M\func{Re}g\left( t\right) \geq \func{Re}f\left( t\right) \geq m\func{Re}%
g\left( t\right)
\end{equation*}%
and%
\begin{equation*}
M\func{Im}g\left( t\right) \geq \func{Im}f\left( t\right) \geq m\func{Im}%
g\left( t\right)
\end{equation*}%
for $\mu $-a.e. $t\in \Omega .$
\end{remark}

If the functions are in $L_{\rho }^{2}\left( \Omega ,\mathbb{R}\right) $
(here $\mathbb{K}=\mathbb{R}$), and $f,g\geq 0,$ $g\left( t\right) \neq 0$
for $\mu $-a.e. $t\in \Omega ,$ then one can state the result:%
\begin{multline}
\int_{\Omega }\rho \left( t\right) f^{2}\left( t\right) d\mu \left( t\right)
\int_{\Omega }\rho \left( t\right) g^{2}\left( t\right) d\mu \left( t\right)
-\left( \int_{\Omega }\rho \left( t\right) f\left( t\right) g\left( t\right)
d\mu \left( t\right) \right) ^{2}  \label{3.11} \\
\leq \left( \sqrt{M}-\sqrt{m}\right) ^{2}\int_{\Omega }\rho \left( t\right)
f\left( t\right) g\left( t\right) d\mu \left( t\right) \int_{\Omega }\rho
\left( t\right) g^{2}\left( t\right) d\mu \left( t\right) ,
\end{multline}%
provided%
\begin{equation}
0\leq m\leq \frac{f\left( t\right) }{g\left( t\right) }\leq M<\infty \quad 
\text{for \ }\mu -\text{a.e. \ }t\in \Omega .  \label{3.12}
\end{equation}

\begin{remark}
\label{r3.4}We notice that (\ref{3.11}) is a generalisation for the abstract
Lebesgue integral of the Klamkin-McLenaghan inequality \cite{KM}%
\begin{equation}
\frac{\sum_{k=1}^{n}w_{k}x_{k}^{2}}{\sum_{k=1}^{n}w_{k}x_{k}y_{k}}-\frac{%
\sum_{k=1}^{n}w_{k}x_{k}y_{k}}{\sum_{k=1}^{n}w_{k}y_{k}^{2}}\leq \left( 
\sqrt{M}-\sqrt{m}\right) ^{2},  \label{3.13}
\end{equation}%
provided the nonnegative real numbers $x_{k},y_{k}$ $\left( k\in \left\{
1,\dots ,n\right\} \right) $ satisfy the assumption%
\begin{equation}
0\leq m\leq \frac{x_{k}}{y_{k}}\leq M<\infty \quad \text{for each \ }k\in
\left\{ 1,\dots ,n\right\}  \label{3.14}
\end{equation}%
and $w_{k}\geq 0,$ $k\in \left\{ 1,\dots ,n\right\} .$

We also remark that Klamkin-McLenaghan inequality (\ref{3.13}) is a
generalisation in its turn of the Shisha-Mond inequality obtained earlier in 
\cite{SM}:%
\begin{equation*}
\frac{\sum_{k=1}^{n}a_{k}^{2}}{\sum_{k=1}^{n}a_{k}b_{k}}-\frac{%
\sum_{k=1}^{n}a_{k}b_{k}}{\sum_{k=1}^{n}b_{k}^{2}}\leq \left( \sqrt{\frac{A}{%
b}}-\sqrt{\frac{a}{B}}\right) ^{2}
\end{equation*}%
provided 
\begin{equation*}
0<a\leq a_{k}\leq A,\text{ \ }0<b\leq b_{k}\leq B
\end{equation*}%
for each $k\in \left\{ 1,\dots ,n\right\} .$
\end{remark}

\end{document}